\documentclass{article}[11pt]
\usepackage{amssymb}
\usepackage{amsmath}
\usepackage{graphicx}
\usepackage{epsfig}
\usepackage[english]{babel}
\title{Statistics of incomplete quotients of continued fractions of quadratic
irrationalities}
\author{E.\,Yu.~Lerner\footnote{Kazan State University, Russia;
e-mail: eduard.lerner@gmail.com}} \textwidth=165mm
\begin{document}
\maketitle
\begin{abstract}
V.I. Arnold has experimentally established that the limit of the
statistics of incomplete quotients of partial continued fractions
of quadratic irrationalities coincides with the Gauss--Kuz'min
statistics. Below we briefly prove this fact for roots of the
equation $r x^2+p x=q$ with fixed $p$ and $r$ ($r>0$), and with
random~$q$, $q\le R$, $R\to \infty$. In Section~3 we estimate the
sum of incomplete quotients of the period. According to the
obtained bound, prior to the passage to the limit, incomplete
quotients in average are logarithmically small. We also upper
estimate the proportion of the ``red'' numbers among those
representable as a sum of two squares.
\end{abstract}


\medskip

\noindent{\bf MSC classes:} 11T06; 11T24; 37E15.

\medskip

\section{The statement of the main results and their consideration}
The papers, monographs, and reports of V.I.~Arnold dedicated to
the statistics of periods of continued fractions and their
multidimensional generalizations (see \cite{ar1}--\cite{probl08})
contain a vast number of experimental facts and outline prospects
for further investigations. The first step in this direction is
the proof of the following assertion (established by V.I.~Arnold):
the statistics of incomplete quotients of the continued fraction of solutions to a random
quadratic equation with integer coefficients, proceeding to the
``thermodynamic'' limit, turns into the Gauss--Kuz'min statistics.
Unfortunately, due to circumstances the proof of this assertion performed by
V.A.~Bykovskii and his followers was not published.
Below we adduce an elementary proof
of the one-parameter version of this assertion.

Any number $x\in {\mathbb R}$ is representable as a
continuous fraction in the form
$$
x=a_0+\cfrac{1}{a_1+\cfrac{1}{a_2+\ldots}} =[a_0;a_1,a_2,\ldots],
$$
where $a_0\in {\mathbb Z}$, $a_i\in {\mathbb N}$ for all $i\ge 1$.

Following Khinchin~\cite{hinchin}, we denote by $E
\binom{s_1,s_2,\ldots,s_k}{A_1,A_2,\ldots,A_k}$ the set of real
numbers which satisfy the conditions
$$a_{s_1}=A_1,\
a_{s_2}=A_2,\ldots,a_{s_k}=A_k;$$
here, certainly, all $A_i$ and
$s_i$ are positive integers, and all $s_i$ are different.
This set is the union of a countable number of intervals~\cite{hinchin}.

Let $I^{s_1,\ldots,s_k}_{A_1,\ldots,A_k}(x)$ stand for the
indicator function of the corresponding set~$E$:
$$
I^{s_1,\ldots,s_k}_{A_1,\ldots,A_k}(x)=\left\{ \begin{array}{ll}
1, & \text{if $x\in E\binom{s_1,\ldots,s_k}{A_1,\ldots,A_k}$,}\\
0, & \text{otherwise,}\end{array}\right.
$$
Let
\begin{equation}
\label{mu} \mu\binom{s_1,\ldots,s_k}{A_1,\ldots,A_k}=\int_0^1
I^{s_1,\ldots,s_k}_{A_1,\ldots,A_k}(x)\,dx
\end{equation}
be the Lebesgue measure of the set of points of the interval $[0,1)$
which belong to $E\binom{s_1,\ldots,s_k}{A_1,\ldots,A_k}$.

Let $x_+(q)$ ``+'' stand for a root of the quadratic equation $r
x^2+p x=q$, i.\,e., $x_+(q)=\frac{-p+\sqrt{\Delta}}{2r}$, where
$\Delta=p^2+4 r q$. In what follows we understand the random
choice as the equiprobable sampling from a finite set.

\medskip
\noindent{\bf Theorem~1.\ \ }{\it Let $r\in{\mathbb N}$,
$p\in{\mathbb Z}$, $q\in \{1,\ldots,R\}$; here $r$ and $p$ are
fixed, $q$ is randomly chosen. Let $P_R^{r,p}
\binom{s_1,\ldots,s_k}{A_1,\ldots,A_k}$ stand for the probability
that $x_+(q)\in E\binom{s_1,\ldots,s_k}{A_1,\ldots,A_k}$. Then
$$
\lim_{R\to\infty}P_R^{r,p}
\binom{s_1,\ldots,s_k}{A_1,\ldots,A_k}=\mu\binom{s_1,\ldots,s_k}{A_1,\ldots,A_k}.
$$
}

\medskip
\noindent{\bf Remark~1.}\ \  {\it One can replace the condition
$q\in \{1,\ldots,R\}$ with that $q\in{\mathbb Z}$, $q\le R$,
$\Delta\ge 0$; evidently, this does not affect the limit
probability.}
\medskip

Using the Kuz'min theorem on the exponential convergence to
$\ln(1+\frac{1}{A(A+2)})/\ln 2$
one can easily obtain the following assertion.

\medskip
\noindent{\bf Corollary of Theorem~1.\ \ }{\it Let $r$, $p$, and $q$
be chosen in the same way as in Theorem~1; let the number $s$, $s\in
\{1,\ldots,n\}$, be chosen randomly. Then as $R\to\infty$ the
probability that $a_s=A$ tends to
$\ln(1+\frac{1}{A(A+2)})/\ln 2+O(1/n)$.}
\medskip

\noindent{\bf Remark~2.\ \  } {\it
The probability $\sum_{s=1}^n P_R^{r,p} \binom{s}{A}/n$
(considered in the corollary) with fixed~$R$ and $n\to\infty$ tends to the probability that
a number chosen randomly from the period of the continued fraction $x_+(q)$ equals $A$.}
\medskip

The proof of Theorem~1 is based on the following fundamental idea:
one can treat integral (\ref{mu}) not only as the Lebesgue integral,
but also as the Riemann integral (considered in the school course on mathematics). A certain special
choice of integral sums actually leads to simplified versions of Theorem~1.
In a general case, it is convenient to perform the proof on the base of the theory of divergent series.

The idea to apply the integral Riemann sums is not new; it was
used implicitly in the proof of a similar correlation in the case
of rational numbers with a fixed denominator~\cite{knuth}.
However, usually one uses another technique for this purpose.
The asymptotic fairness of the Gauss--Kuz'min statistics in the case of
a fixed denominator follows from results obtained by L.~Heilbronn and J.W.~Porter
(see \cite{UstinovAlgebraAnaliz}). Not long time ago,
V.A.~Bykovskii, M.O.~Avdeeva, A.V.~Ustinova \cite{bik}, \cite{avd}, \cite{us}
obtained the limit statistics for finite fractions, whose numerators and
denominators belong to a sector or to an arbitrary expanding domain.

Let us explain the (mentioned in Remark~2) connection with the statistics of the period
of a continued fraction. In accordance with the Lagrange theorem, a continued fraction for the quadratic
irrationality (and only for it) is periodic. Therefore the continued
fraction for $x_+(q)$ takes the form
\begin{eqnarray}
\frac{-p+\sqrt{p^2+4rq}}{2r}&=&[a_0;\ldots,a_m,a_{m+1},\ldots,a_{m+T},a_{m+1},\ldots,a_{m+T},\ldots]=\nonumber\\
\label{CFforQI} &=&[a_0;\ldots,a_m,[a_{m+1},\ldots,a_{m+T}]],
\end{eqnarray}
where $T=T(r,p,q)$ is the length of the period of the continued
fraction, $m=m(r,p,q)$ is the length of the preperiod. Evidently,
one can represent the probability that with random $q$, $q\le R$,
a number randomly chosen from the period $a_{m+1},\ldots,a_{m+T}$
coincides with~$A$ as the limit $\lim\limits_{n\to
\infty}\sum_{s=1}^n P_R^{r,p} \binom{s}{A}/n$. Unfortunately, we
did not succeed in proving that this probability tends to
$\ln(1+\frac{1}{A(A+2)})/\ln 2$ as $R\to\infty$.

In papers \cite{umn},\cite{ivestiya} V.I.~Arnold experimentally
studies~$T(1,p,q)$. He establish that in average this period
increases proportionally to the square root of the discriminant:
\begin{equation}
\label{simArnold} \widehat{T} \sim \text{const} \sqrt{\Delta}.
\end{equation}

Let $T_0(q)=T(1,0,q)$ be the period of the continued fraction of
the square root of~$q$. Let $\widehat{T}'(Q)=\sum_{q=1}^Q
T_0(q)/Q$. This kind of an average was studied experimentally by experts in
the theory of numbers, who stated the conjecture
\begin{equation}
\label{down} \widehat{T}_0(Q) \sim \text{const} \sqrt{Q}
\log^\alpha (Q),\quad \text{where }\alpha <0
\end{equation}
(see \cite{Gus12} and references therein).

One can easily obtain the following upper bound for $\widehat{T}_0(Q)$ ~\cite{av}:
$$\widehat{T}'(Q) <
\text{const} \sqrt{Q}.
$$
In \cite{GusMS} (see also \cite{Gus12}) one proves that
the left-hand side of this inequality is asymptotically small in comparison with the right-hand
one, moreover, if the extended Riemann hypothesis is fulfilled, then the order
of their difference is not less than $\log(Q)^{\log
2-\varepsilon}$.

The well-known lower estimate for $\widehat{T}_0(Q)$ significantly
differs from the right-hand side of~(\ref{down}). In \cite{Gus12}
one only proves that $\widehat{T}_0(Q) > \text{const} \log(Q)$. Note
that obtaining the bound $\widehat{T}_0(Q)> \text{const} \sqrt{Q}
\log^\alpha (Q)$ would have allowed one to solve the Gauss
question~\cite{venkov} about the slow growth of the number of
classes of real quadratic fields.

There is a known bound for the maximal value of $T_0(q)$. Put
\begin{equation}
\label{DQ}
D(Q)=\sum_{u=1}^{\lfloor\sqrt{Q}\rfloor}\tau(Q-u^2),
\end{equation}
where $\lfloor\cdot\rfloor$ is the integer part, $\tau(n)$ is the
quantity of divisors of the number~$n$. Formally speaking, the
Dirichlet theorem
(the equality $\sum_{u=1}^Q \tau(u)/Q=\ln
Q+(2\gamma-1)+O(1/\sqrt{Q})$) {\bf does not imply} that $D(Q)\sim
\sqrt{Q} \ln Q$, because the values of the function $\tau$ are not uniform.
Reasoning more accurately (see \cite{av} and references therein), we obtain the bound
\begin{equation}
\label{D} D(Q)=O(\ln^3(Q)\sqrt{Q}).
\end{equation}

The mentioned bound for the maximal value of $T_0(q)$ (recently it
was described by V.I.~Arnold in~\cite{vid3}) has the form
$$
T_0(Q)\le D(Q).
$$
This inequality was first obtained by Hickerson in
1973~\cite{hickerson} (see also \cite{pen}).
Later Western researchers
~\cite{cohn}, \cite{stanton} succeeded in proving the correlation
\begin{equation}
\label{piter} T_0(Q)=O(\sqrt{Q} \ln Q).
\end{equation}
Earlier mathematicians of the Leningrad scientific school treated it
as a direct corollary of the Dirichlet formula and results of~\cite{pen} (see
\cite{lomi}). Moreover, if the extended Riemann hypothesis is fulfilled, then one can analogously obtain the
double logarithmic bound
$$
T(r,p,q)=O(\sqrt{\Delta} \ln\ln \Delta).
$$

One can construct a continued fraction for a quadratic irrationality,
using the algorithm of successive approximations by
mediants (the Farey fractions). Based on the connection
of this algorithm with explicit representations of values of an integer-valued quadratic
form (see \cite{conway}), in Section~3 we obtain the following
bound for the sum of partial quotients of the quadratic irrationality within
one its period.

\medskip

\noindent{\bf Theorem~2.}\ \ {\it Let $r x^2+p x=q$ be an
arbitrary quadratic equation with integer coefficients and
real irrational roots. In denotations of
(\ref{CFforQI}), (\ref{DQ}) we have the inequality
\begin{eqnarray}
\label{th2} \sum_{i=1}^{T(r,p,q)} a_{m(r,p,q)+i}\le
f(\Delta/4),\qquad
\text{ where\phantom{nbbnbnbnb}}\\
\label{th2even} f(Q)= 2 D(Q)+\tau(Q), \text{if $Q$ is integer;}\\
\label{th2odd}
f(Q)=2\mathop{\sum\limits_{i\in\{1,3,\ldots\}}}\limits_{i^2<4Q}\tau(Q-i^2/4),
\text{  otherwise.   }
\end{eqnarray}
In the case of an odd period~$T(r,p,q)$ one can improve this
bound, namely, divide the right-hand side of
inequality~(\ref{th2}) by two.}
\medskip

According to the results of computer experiments, for almost all roots of
prime numbers in the form $4k+3$, as well as the roots of the doubled numbers
in the same form, the left-hand side of inequality~(\ref{th2}) is
half as large as the right-hand one.

Problem~5 in the list of problems stated by Arnold
in~\cite{probl08} implies the estimation of rate of growth of the
average value of elements of the period of a continued fraction.
To put it more precisely, the problem is to solve the alternative
between the power rate of growth and the logarithmic one.
Comparing inequality~(\ref{th2}) with formula~(\ref{D}) and the
results of the numerical tests (\ref{simArnold},\ref{down}), we
conclude that the mean value $\widehat a$ is logarithmically small
as against~$\Delta$.

Note that the geometric method used for the proof of Theorem~2
enables one to easily obtain assertions on the palindromicity of
continued fractions in cases when $p=0$ or $r=1$ (see
remarks~5~and~6). We are going to dedicate a separate paper to the study of problem~12
in~\cite{probl08}.


Studying periods lengths $T_0(Q)$, V.I.~Arnold stated the problem
of the constructive definition of the set~$K$ of the so-called
``red'' numbers:
$$
K=\{ Q: Q\in{\mathbb N}, T_0(Q) \text{\ is odd}\}.
$$
The set~$K$ is a subset of the totality~$M$ of positive integers
representable as a sum of two squares. The explicit representation
of the set~$M$ is well known (see, for example, \cite{davenport}). In
particular, prime numbers $Q$ of the set~$M$ obey the formula:
$Q\!\!\mod\! 4\,\ne 3$.

Let~$K_n$ (or~$M_n$) stand for the part of the set~$K$
(the set~$M$) which consists of numbers not greater than~$n$. According to the
Arnold conjecture, the following limit exists:
$$
\lim_{n\to\infty} \frac{K_n}{M_n}=c.
$$

The sets~$K$ and $M$ are closely connected with the generalized Pell equation
\begin{equation}
\label{noPell} x^2-Q y^2=-1.
\end{equation}

\noindent{\bf Theorem~3.}\newline {\bf 1.}\ \ {\it $Q\in K
\Leftrightarrow$ equation~(\ref{noPell}) is solvable in integer
numbers.}
\newline {\bf 2.}\ \ {\it $Q\in M
\Leftrightarrow$ equation~(\ref{noPell}) is solvable in rational
numbers.}
\newline {\bf 3.}\ \ {\it Let $P$ be the set of prime numbers. With sufficiently
large~$n$,}
$$|K_n|<|M_n| \mathop{\prod_{p:p\in P,}}_{p \!\!\!\!\mod\! 4\,\ne 1}(1-\frac{1}{p^2}).$$
\newline {\bf 4.}\ \
{\it $Q\in K\cap P\quad\Leftrightarrow\quad Q\in P,\ \ Q\!\!\mod\! 4\,\ne 3\quad\Leftrightarrow\quad Q\in M\cap P$.}

\medskip

One can easily calculate that
$$
\mathop{\prod_{p:p\in P,}}_{p \!\!\!\!\mod\! 4\,\ne
1}(1-\frac{1}{p^2})=0.64208\ldots.
$$
This gives a rough bound for~$c$; computer tests performed for all
$n$ from $10^5$ to $10^6$ show that $0.47<\frac{K_n}{M_n}<0.48$.

Theorem~3.1 is well known for experts in the theory of numbers
(see, for example, \cite[Theorem 12]{venkov}), as well as the absence
of simple solvability criteria for
equation~(\ref{noPell}) (i.\,e., for the membership of~$K$).


The assertion about the sufficiency in Theorem~3.2 is evident. The
proof of the necessity is reduced to the proof of the following
property: if a positive integer is representable as a sum of two
squares of rational numbers, then it is representable as a sum of
two squares of integers. The proof for the case of two squares
verbatim coincides with that adduced in~\cite[Appendix to
Chapter~4]{serr} for the case of three squares.

Theorem~3.3 follows from two evident corollaries of the previous items of
Theorem~3.  Apparently, Theorem~3.2 implies that the set~$M$ contains
sets of all numbers in the form $4 K$, $9 K$, $49 K$, etc
(here we restrict ourselves only with the multiplication by squares
of prime numbers, for which the residue of division by~4 differs from~1).
On the other hand, Theorem~3.1 implies that $K$ does not intersect with these sets.
To put it more precisely, if $Q\in K$, then $Q$ is neither divisible by 4 nor
by all prime $p$ in the form $4k+3$. Otherwise, considering
equality~(\ref{noPell}) modulo~4 or modulo~$p$, we get the
equation $x^2\equiv -1$; it is well known that the latter equation is
neither solvable modulo~4 nor modulo~$p$.

Now Theorem~3.4 is reduced to the assertion that all primes in the
form $4k+1$ belong to~$K$. This fact was proved by Legendre~(see
\cite{davenport}, See~\cite{venkov}) for the proof.
He used it in order to solve the inverse problem
on the representation of such a number as a sum of two squares.

Note that Theorem~3 is absent in author's article with the same
title submitted to the journal ``Functional Analysis and Its
Applications''

This paper is written under the bright impression of (video)
lectures delivered by V.I.~Arnold
\cite{vid1},~\cite{vid2},~\cite{vid3}. It resulted from the
addition to these lectures which, in turn, appeared as a result of
the performed numerical experiments and the study of remarkable
books \cite{knuth}, \cite{conway}, \cite{venkov}.

\section{How to apply the condition ``accurate to a zero-measure set''
in a discrete case (the proof of Theorem~1)}

Let us first prove Theorem~1 in the case of the simplest quadratic
irrationality~$\sqrt{q}$.

\medskip

\noindent{\bf Lemma 1.}\ \ {\it Let $q$ be arbitrarily
chosen from the set $\{n^2+1,n^2+2,\ldots,n^2+2n\}$, $n\in {\mathbb N}$.
Denote by $P'_n \binom{s_1,\ldots,s_k}{A_1,\ldots,A_k}$ the
probability that $\sqrt{q}\in
E\binom{s_1,\ldots,s_k}{A_1,\ldots,A_k}$. Then
$$
\lim_{n\to\infty}P'_n
\binom{s_1,\ldots,s_k}{A_1,\ldots,A_k}=\mu\binom{s_1,\ldots,s_k}{A_1,\ldots,A_k}.
$$
}

{\bf Proof of Lemma 1.} Let us first note that under the assumptions of the
lemma $\lfloor q \rfloor=n$. Therefore, by definition we have
\begin{equation}
\label{def}
P'_n\binom{s_1,\ldots,s_k}{A_1,\ldots,A_k}=\frac{1}{2n}\sum_{i=1}^{2n}
I^{s_1,\ldots,s_k}_{A_1,\ldots,A_k}(\sqrt{n^2+i}-n).
\end{equation}

The complement of the set $E \binom{s_1,\ldots,s_k}{A_1,\ldots,A_k}$
is representable in the form $ \bigcup\limits_{(B_1,\ldots,B_k):
B_j\ne A_j}E \binom{s_1,\ldots,s_k}{B_1,\ldots,B_k} $, i.\,e., as the
union of a countable number of intervals. Consequently, the function
$I^{s_1,\ldots,s_k}_{A_1,\ldots,A_k}$ is countably continuous, therefore,
by the Lebesgue theorem one can treat integral~(\ref{mu})
not only as the Lebesgue integral, but also as the Riemann integral.
In order to calculate it, we divide the segment $[0,1)$ onto $2n$ equal parts.
We choose $\sqrt{n^2+i}-n$ as the point, where we calculate the integrand on the
interval $\left[\frac{i-1}{2n},\frac{i}{2n}\right)$
(one can easily make sure that $\frac{i-1}{2n}<\sqrt{n^2+i}-n<\frac{i}{2n}$ for all
$i=1,\ldots,2n$). Then the value of the Riemann integral sum
coincides with the right-hand side of equality~(\ref{def}). Lemma~1 is proved.

Now the proof of Theorem~1 in the case $r=1$, $p=0$ follows from several remarks.
\medskip

\noindent{\bf Remark~3.}\ {\it For an arbitrary sequence $s_n$ and
any sequence of nonnegative numbers $t_n$ such that
$\sum_{n=1}^{\infty}t_n=\infty$, the convergence $s_n\to\mu$
implies that
\begin{equation}
\label{R}
 S_n=\frac{t_1 s_1+\ldots+t_n s_n}{t_1+\ldots+t_n}\to\mu.
\end{equation}
In the theory of divergent series this property is known as the
regularity of the Riesz summation method~\cite[Chap.~3,
Theorem~12]{hardy}.}

\medskip

\noindent{\bf Remark~4.}\ {\it The assertion of Theorem~1 on the
convergence of a sequence~$P_n$ to $\mu$ follows from the
convergence of a subsequence $P_{n_i}$ to $\mu$, if
$\lim\limits_{i\to\infty} \frac{n_i}{n_{i+1}}\to 1$.}

\medskip

Really, it is evident that
$$\frac{n_i}{R} P_{n_i}\le P_R\le \frac{n_{i+1}}{R} P_{n_{i+1}}\quad
\text{with } R\in \{n_i,n_i+1,\ldots,n_{i+1}\}.$$ In Remark~4 and
in what follows we omit fixed parameters of functions~$P$, $\mu$,
and~$P'$.

In order to complete the proof of Theorem~1 in the case $r=1$, $p=0$, it remains
to note that by definition
 $P_{R^2+2R}^{1,0}=\sum_{n=1}^R
\frac{2n}{R^2+2R} P'_n$, i.\,e.,
\begin{equation}
\label{temp} \frac{R^2+R}{R^2+2R}P_{R^2+2R}^{1,0}=\sum_{n=1}^R
\frac{2n}{R^2+R} P'_n
\end{equation}
According to Lemma~1 and Remark~3, the limit of the right-hand
side of equality~(\ref{temp}) for $R\to\infty$ coincides with the
corresponding value of $\mu$, and the limit of the left-hand side
(due to Remark~4) is $\lim\limits_{R\to\infty} P_R^{1,0}$.

One can prove Theorem~1 analogously in the case $r=p=1$. According
to~\cite{ivestiya}, one can reduce the general proof with~$r=1$ to these two cases.
Unfortunately, partitions of the integration interval $[0,1)$ onto
equal parts do not enable us to prove Theorem~1 with $r>1$. More intricate techniques are necessary.

{\bf Proof of Theorem~1 in the general case.} Let $n\in{\mathbb
N}$. In order to calculate integral (\ref{mu}), let us divide the
segment $[0,1)$ onto $p+r+2nr$ parts with left ends at the points
$x_+(q)-n$; here $q=i-1+n(p+n r)$, $i$ stands for the number of an
interval. Let us calculate the value of the integrand on the
interval at the same point~$x_+(q)-n$ (further $x(q)\equiv
x_+(q)$). We obtain
\begin{equation}
\label{temp2} \mu=\lim_{n\to\infty} \sum_{q=n(p+n
r)}^{(n+1)(p+(n+1) r)-1} (x(q+1)-x(q))\ I(x(q)).
\end{equation}
In accordance with Remark~3 and due to the boundedness of the sums
in the right-hand side of~(\ref{temp2}), we have
$$\mu=\lim_{n\to\infty} \sum_{q=1}^{n(p+n r)}
\frac{x(q)-x(q-1)}{n} I(x(q-1))=\lim_{N\to\infty} \sum_{q=0}^{N-1}
\frac{x(q+1)-x(q)}{\lfloor x(N) \rfloor} I(x(q)).
$$
Evidently, in the denominator of the latter expression we can write $x(N)$
in place of $\lfloor x(N) \rfloor$. As a result, we get the limit
in form (\ref{R}), where $t_i=x(i)-x(i-1)$. It is well known that (see
\cite[Chap.~3, theorem~14]{hardy})
$$
 S_N=\frac{t_1 I_1+\ldots+t_N I_N}{t_1+\ldots+t_N}\to\mu\quad
\Leftrightarrow\quad
  S'_N=\frac{t'_1 I_1+\ldots+t'_N I_N}{t_1+\ldots+t_N}\to\mu,
$$
provided that the following conditions hold:
\begin{eqnarray}
\nonumber
\frac{t_{N+1}}{t_N}\le \frac{t'_{N+1}}{t'_N},\\
\nonumber \sum_{q=1}^N t_q/t_N\le \text{const} \sum_{q=1}^N
t'_q/t'_N.
\end{eqnarray}
Let us verify them in the case $t_q=x(q)-x(q-1)$, $t'_q\equiv
1$. The first condition $x(q-1)+x(q+1)\le 2x(q)$
follows from the convexity of the function $x(q)$,
the second one does from the existence of the limit
$\lim\limits_{N\to\infty} N \left(
1-\frac{x(N-1)}{x(N)}\right)=1/2$. Therefore, we get
$$
\mu=\lim_{N\to\infty} \sum_{q=0}^{N-1} I(x(q))/N=\lim_{R\to\infty}
\sum_{q=1}^{R} I(x(q))/R.
$$
By definition $\sum_{q=1}^{R} I(x(q))/R=P_R^{r,p}$. Theorem~1
is proved.

\section{The gradual ``Nose-Hoover'' algorithm (the proof of Theorem~2)}
{\it \small ~\quad\quad\quad\quad``\ldots
 and before he knew what he was doing he lifted up his trunk and hit that fly dead with the end of it.''\\
\phantom{mmmmmmmmmmmmmmmmmmmmmmmmmmmmmmmmmmmmmmmmmmm} R.Kipling.
``The Elephant's Child'' }

\medskip
In what follows, without loss of generality, we consider only continued fractions of positive irrational numbers.

The standard ``Nose-Hoover'' algorithm for finding a continued
fraction of a real number $x$ was proposed by Klein in
1895~\cite{davenport} (the term was introduced by B.N.~Delone). It
is well known~\cite{ar1} that this algorithm is reduced to the
geometric method which constructs the boundary of the convex shell
of the set of integer nonnegative points~$(u,v)$ located above
(below) the straight line~$v=x u$. Let $\pmb{e}_0$ stand for the
vector $(0,1)$, $\pmb{e}_1=(1,0)$. Put
\begin{equation}
\label{step} \pmb{e}_{n+1}=\pmb{e}_{n-1}+a_{n-1} \pmb{e}_n,\quad
n=1,2,\ldots,
\end{equation}
where $a_{n-1}$ is the maximal integer such that the vector
$\pmb{e}_{n+1}$ lies below (for odd $n$) or above (for even $n$)
the straight line~$v=x u$. Klein noted that $a_n$, $n=0,1,\ldots$,
coincide with partial quotients of the continued fraction of the
number~$x$. The geometric algorithm results in two polylines which
represent parts of sails of the continued fraction.

We need a slightly modified algorithm which results in one
infinite polyline~$L$, originating at the point~$(1,1)$, whose
segments are the vectors $a_{n-1} \pmb{e}_{n}$, $n=1,2,\ldots$.
Evidently, $\pmb{e}_{n+1}+\pmb{e}_n=\pmb{e}_n+\pmb{e}_{n+1}$.
Consequently, if in the standard ``Nose-Hoover'' algorithm
~(\ref{step}) we add the extra vector $\pmb{e}_n$ and thus go out
of the line $v=x u$ (i.\,e., we consider the vector
$\pmb{e}_{n-1}+(a_{n-1}+1) \pmb{e}_n$), then we get the first
integer point on the segment of the polyline of another sail
constructed at step $n+1$. This fact justifies a simple geometric
algorithm which constructs the polyline~$L$.

We begin the construction process with the point $(1,1)$; it is
convenient to connect it with the origin of coordinates by the
segment which does not enter in~$L$. As the  main ``constructive''
term at the first ($n$th) step we choose the vector $\pmb{e}_1$
(the vector $\pmb{e}_n$). We add this vector till we go out of the
line $v=x u$. The newly added vector $\pmb{e}_{n+1}$ is directed
from the origin of coordinates to the point obtained as a result
of the latter addition up to the step out of the line. See Fig.~1
for the first segments of the polyline for $x=\sqrt{2}$.

For convenience, we assume (if the contrary is not specified)
that the zero partial quotient (i.\,e., $\lfloor x
\rfloor$) differs from zero. In order to provide the correspondence to
the indices of partial quotients, it is also convenient to begin the numeration
of segments of the polyline~$L$ with zero.

\begin{figure}
\begin{center}
\includegraphics{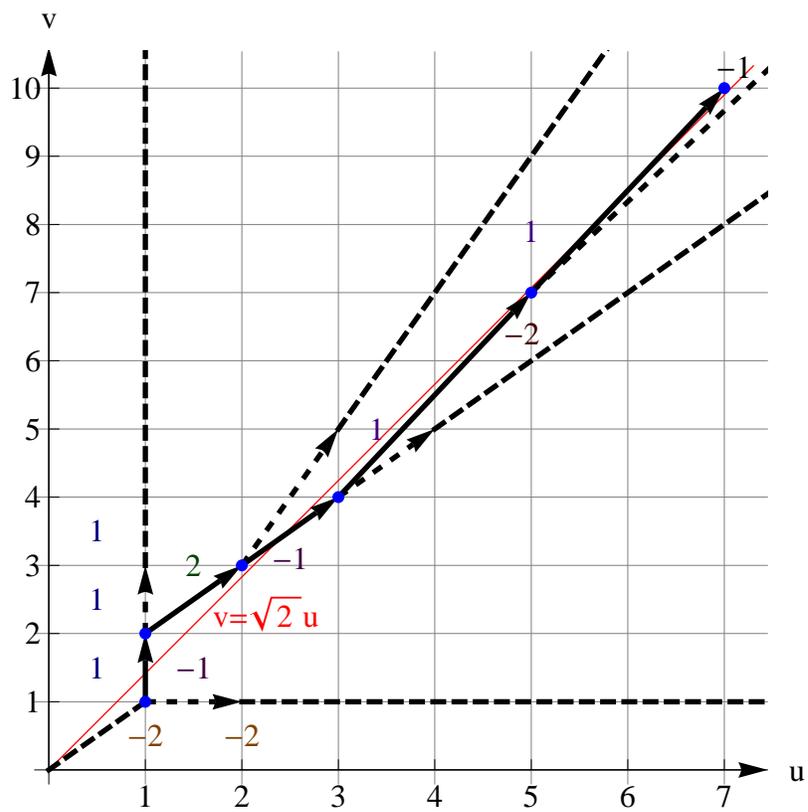}
\end{center}
\caption{The gradual ``Nose-Hoover'' algorithm `` for $x=\sqrt{2}$} \label{1}
\end{figure}

Let us first consider the approximation of a number~$x$ by
mediants. The study of series of mediants is usually connected
with the names of the German mathematician Stern, the French
watchmaker Brocot (he described this series and applied it in
clock manufacturing) or (more often) the English geologist Farey.
The latter paid attention to the ``curious property of vulgar
fractions'' (this was the title of a short letter of Farey
published in ``Philosophical Magazine'' in 1816). Actually, the
Farey sequence was described by Haros (see
\cite[p.~36--37]{hardywright}); he considered its application in
1802. Even ancient Greek mathematicians were able to enumerate
rational numbers with the help of mediants.

Let $f_1=\frac{v_1}{u_1}$, $f_2=\frac{v_2}{u_2}$ be irreducible
fractions such that $v_1,u_1,v_2,u_2\ge 0$ and $f_1<f_2$. Let us
define the operation $\downarrow$ of finding the mediant (the
``insertion'' operation) by the formula $f_1\downarrow
f_2=\frac{v}{u}$, where $v=v_1+v_2$, $u=u_1+u_2$. Earlier we used
the denotation introduced by A.A.~Kirillov. Evidently, in the
geometric representation of the fraction $\frac{v}u$ as the
integer vector $(u,v)$ the operation $\downarrow$ corresponds to
the addition of vectors $(u_1,v_1)$ and $(u_2,v_2)$. The obtained
diagonal of the parallelogram appears to be ``inserted'' between
its sides.

Let $x\in (f_1,f_2)$, $f_3=f_1\downarrow f_2$. Since $f_3\in
(f_1,f_2)$, one can treat $f_3$ as an approximation of the
number~$x$. If $x\in (f_1,f_3)$, then we put $f_4=f_1\downarrow
f_3$, otherwise we do $f_4=f_3\downarrow f_2$. The process of the
approximation of the number~$x$ by mediants consists in the
repetition of these operations for the corresponding intervals.

The condition $x>0$ means that $x\in\left( \frac01,\frac10 \right)$.
The initial approximation of $x$ for such an interval is the
fraction~$\frac11$. One can easily see that for the mapping
$$ \text{the fraction $\frac{v}u \quad \leftrightarrow\quad$ the point
$(u,v)$}
$$
the geometric representation of the algorithm for approximating
the number~$x$ by mediants is completely identical to the
algorithm for constructing the polyline~$L$. The technical
distinction consists in the following fact: earlier we constructed
each segment of the polyline~$L$ ``at once'', but now it ``grows
gradually'' due to the stepwise addition of the next
vector~$\pmb{e}_i$. We call a simplified version of this algorithm
as applied to quadratic irrationalities the gradual
``Nose-Hoover'' algorithm.

Let $(u_n,v_n)$ stand for the $n$th integer point on the
polyline~$L$, $n=0,1,2,\ldots$. We obtain it at the $n$th step of
the gradual ``Nose-Hoover'' algorithm. In Fig.~1, for example,
$(u_0,v_0)=(1,1)$, $(u_1,v_1)=(1,2)$, $(u_2,v_2)=(2,3)$,
$(u_3,v_3)=(3,4)$, $(u_4,v_4)=(5,7)$, etc. The number~$i$ of the
vector $\pmb{e}_i$ which we add at the $n$th step, generally
speaking, does not exceed~$n$; the identity $i\equiv n$ takes
place only for the golden ratio.

The main inconvenience of the standard ``Nose-Hoover'' algorithm consists
in the following fact: the polylines very quickly tend to the straight line $u=x v$
and the geometric illustration requires the high accuracy (``all noses
become longer''). In the case of a quadratic
irrationality we can use only conventional illustrations which allow us to
``visualize'' the topology of polylines~$L$ (and therefore the continued
fraction) without its exact representation on the integer lattice.
In essence, these illustrations coincide with the approach (described in~\cite{conway})
to the visualization of values of a quadratic form on the banks of its river
(we use the terminology of~\cite{conway}). Below we show the geometric juxtaposition.

Let us draw one more vector which originates from the point
$(u_n,v_n)$; denote it by $\pmb{e}'_n$. The vector $\pmb{e}'_n$ is
defined by the condition $\pmb{e}'_n+\pmb{e}_i+\pmb{e}''_n=0$,
where $\pmb{e}_i$ is the vector which originates from the point
$(u_n,v_n)$ and goes along the polyline~$L$,
$\pmb{e}''_n=-(u_n,v_n)$ (we considered this vector earlier). The
collection of vectors $\{\pmb{e}'_n,\pmb{e}_i,\pmb{e}''_n\}$,
where each one is defined accurate to the multiplier $(-1)$, is a
{\it superbasis}~\cite{conway}. This means that any pair of these
vectors generates the whole integer lattice. The transition from
the point $(u_n,v_n)$ to that $(u_{n+1},v_{n+1})$ corresponds to
the replacement of one superbasis with another one; the latter
differs from the initial superbasis only in one of three its
elements. Thus, for the polyline~$L$ in Fig.~1 the initial
superbasis is $\{(0,1),(1,0),(1,1)\}$, the 1st one is
$\{(1,1),(1,0),(1,2)\}$, the 2nd one is $\{(1,2),(1,1),(2,3)\}$,
the 3rd one is $\{(1,1),(2,3),(3,4)\}$, etc.

Let us extend vectors~$\pmb{e}'_n$ up to the rays which originate
at the points $(u_n,v_n)$. As a result, the first quadrant appears
to be divided onto disjoint connected domains (see Fig.~1). Let us
associate each vector $(u_n,v_n)$ with a domain, whose boundary
contains all points of polylines which include the vector
$(u_n,v_n)$ in the corresponding superbasis. In the picture these
domains look like narrow ``crevices'' located at the north-east of
the points $(u_n,v_n)$. We associate the domains which border on
the coordinate axes with the unit vectors of these axes.

Without loss of generality, in what follows we assume that in this
statement of Theorem~2 the coefficient $r$ is positive.

If, in addition, the equation has a unique positive root, then
\begin{equation}
\label{equiv} \text{$the point (u_n,v_n)$ is located above the line
$v=x_+(q)u$}\ \Longleftrightarrow\ r v_n^2+p v_n u_n-q u_n^2>0.
\end{equation}
If both roots are positive, then this equivalence does not
necessarily take place for all~$n$. The quadratic form is also
positive at the points located below the line $v=x_-(q) u$.
However, since the polyline~$L$ tends to the straight line
$v=x_+(q) u$ (it is known that $|\frac{v_n}{u_n}-x|<
1/u_{n-1}^2$), correlation (\ref{equiv}) is true for all $n\ge
n_0$. Let $l_0$ stand for the number of the segment of the
polyline which contains the point $(u_{n_0},v_{n_0})$. In the case
of one positive root, $l_0=0$.

Let us write the values of the quadratic form $r v^2+p v u-q u^2$ at the points
$(u_n,v_n)$ in the corresponding domain. See Fig.~1 for the result obtained
in the case $r=1$, $p=0$, $q=2$ (the upper right corner of the figure contains
the value calculated at the point~$(5,7)$). From~(\ref{equiv}) it follows
that with all $n\ge n_0$ the points of the polyline~$L$ belong to boundaries of
domains with positive and negative values. The collection of
superbases which correspond to these points is called the {\it river} of the
quadratic form~\cite{conway}.

Let us adduce the main result of the theory of quadratic forms which enables
one to make the calculation of the period of a continued fraction much easier.
Using the values of the quadratic form at three vectors from any superbasis,
one can easily restore the continued fraction. To put it more precisely, the following
correlation~\cite{conway} is true. Let domains associated with values
$a,b,c,d$ of the quadratic form be located in accordance with Fig.~2.
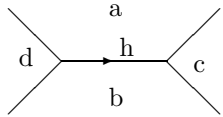
\begin{figure}[t]
\begin{center}
\begin{picture}(100,40)
\put(30,20){\line(-1,1){20}} \put(30,20){\line(-1,-1){20}}
\put(30,20){\vector(1,0){20}} \put(50,20){\line(1,0){20}}
\put(70,20){\line(1,1){20}} \put(70,20){\line(1,-1){20}}
\put(14,18){d} \put(80,16){c} \put(48,37){a} \put(48,3){b}
\put(52,22){h}
\end{picture}
\end{center}
\caption{The arithmetic progression rule} \label{1}
\end{figure}
Then values $d,a+b,c$ form an arithmetic progression.
Moreover, if $h$ is the common difference of this progression, then
\begin{equation}
\label{main} h^2/4-a b= \Delta/4.
\end{equation}

Therefore, using the values $(a,b)$ for two neighboring domains
with different signs separated by the polyline~$L$, and the value
$h$, one can unambiguously restore all subsequent values of
partial quotients of the continued fraction, moving along the
river of the quadratic form (see the figures in Chapter 1 of book
\cite{conway}). Moreover, according to the arithmetic progression
rule, one can restore all previous values of partial quotients
(whose numbers exceed~$l_0$).

\medskip
\noindent{\bf Remark~5.}{\ \it If $h=0$, then the subsequent
values and the previous ones coincide. In this case the periodic
part of the continued fraction is palindromic. In particular, this
takes place for $p=0$, i.\,e., for the root of a rational number
(see the latter paragraph of this section).}

\medskip
\noindent{\bf Remark~6.}{\ \it The uniqueness of the river of a
quadratic form~(see~\cite{conway}) implies that indefinitely
restoring the previous values, we obtain (beginning with certain
$n'_0$) the basic (infinite) part of the continued fraction of the
associated irrationality~$x_-(q)$. Hence, one can obtain the
period of the continued fraction~$x_-(q)$ from the period
$a_{m(r,p,q)+1},\ldots,a_{m(r,p,q)+T}$ by operations of the
inversion and the cyclic shift. In particular, for $r=1$, when
$x_-(q)$ is connected with $x_+(q)$ by the evident unimodular
transform (it is well known that the latter preserves the period
of a continued fraction), we observe the palindromicity.}
\medskip

\begin{table}[t]
\begin{center}
\begin{tabular}{|c|c|c|c|c|c|}
\hline
{\it the number of a segment of the polyline } &\it0&\it1&\it1&\it2&\it2\\
\hline
$a$ is the value located at the north west of~$L$&1&2& 1&1&1\\
\hline
$b$ is the value located at the south east of~$L$&-1&-1&-1&-2&-1\\
\hline
$h$ is the common difference of the arithmetic progresion&2&0&-2&0&2\\
\hline
\end{tabular}
\end{center}
\caption{All possible values of $(a,b,h)$ for the polyline
represented in Fig.~1}
\end{table}


Let us now immediately prove the theorem. Its
assertion directly follows from equality~(\ref{main}).
Really, as was mentioned above, one can perfectly restore partial
quotients~$a_i$ with $i>l_0$, using the gradual
``Nose-Hoover'' algorithm with $n\ge n_0$. This can be done not only with the help
of the exact geometric representation in the coordinate plane, but also
by a ``caricature'' representation of the motion along the river of the quadratic form.
Assume that at step~$n_1$ we get the same parameters of the arithmetic progression
as those obtained at step~$n_0$. Let~$l_1$ stand for the number of the polyline segment
which contains the point $(u_{n_1},v_{n_1})$. Then
a part of the continued fraction $a_{l_0+1},a_{l_0+2},\ldots,a_{l_1}$ becomes
periodic. In addition, the number $l_1-l_0$ is even, because
the points $(u_{n_0},v_{n_0})$ and $(u_{n_1},v_{n_1})$ are located
to one side of the polyline~$L$. In the case of an odd period the
sequence $a_{l_0+1},a_{l_0+2},\ldots,a_{l_1}$ contains
at least two repeating subsequences. In accordance with the
algorithm we have $\sum_{i=l_0+1}^{l_1} a_{i}=n_1-n_0$. Therefore,
suffice it to estimate the number of all possible triplets $(a,b,h)$ which
satisfy~(\ref{main}). The obtained bound
obeys formulas~(\ref{th2even}), (\ref{th2odd}). The multiplier~2
in the right-hand sides of these formulas appears, because we have to take into account various
signs of~$h$, and the term
$\tau(\Delta/4)$ does when we consider the case
$h=0$. Theorem~2 is proved.

Let us consider the case of one positive (and one
negative) root, which is equivalent to $r q>0$. If $\lfloor
x_+(q) \rfloor=0$, i.\,e., $a_0=0$, then the numeration of segments of the
polyline shifts by one. The above considerations imply that only
three alternatives are possible. Namely, a continued fraction can have no
preperiod ($m(r,p,q)=0$), or $m(r,p,q)=1$, or $m(r,p,q)=2$
and $a_0=0$. In the second and third cases the latter partial quotient of the
preperiod is less than the latter partial period of the continued
fraction.


\end{document}